\documentclass[pdflatex]{sn-jnl}

\usepackage{graphicx}%
\usepackage{amsmath,amssymb,amsfonts}%
\usepackage{xcolor}%
\usepackage{parskip}
\usepackage[numbers,sort&compress]{natbib}
\hypersetup{
  breaklinks=true,
  colorlinks=true,
  linkcolor=blue,
  citecolor=blue,
  filecolor=blue,
  urlcolor=blue,
  menucolor=blue
}
\usepackage{enumitem}

\theoremstyle{thmstyleone}%
\newtheorem{theorem}{Theorem}[section]
\newtheorem{remark}{Remark}%

\theoremstyle{thmstylethree}%
\newtheorem{definition}{Definition}[section]%

\newtheorem{lemma}[theorem]{Lemma}
\numberwithin{equation}{section}

\begin{document}

\title[Relations in the Hybrid Kuramoto Flow]{{\color{black} Synchronization Relations in the Hybrid Kuramoto Flow: Equivalence and a High-Coherence Criterion}}


\author*[1]{\fnm{Ting-Yang} \sur{Hsiao}}\email{thsiao@sissa.it}

\author[2]{\fnm{Yun-Feng} \sur{Lo}}\email{ylo49@gatech.edu}

\author[3]{\fnm{Chengbin} \sur{Zhu}}\email{cz43@illinois.edu}

\affil*[1]{\orgdiv{Mathematics}, \orgname{International School for Advanced Studies (SISSA)}, \orgaddress{\street{via Bonomea 265}, \postcode{34136}, \state{Trieste}, \country{Italy}}}

\affil[2]{\orgdiv{School of Electrical and Computer Engineering}, \orgname{Georgia Institute of
Technology}, \orgaddress{\street{777 Atlantic Drive NW}, \city{Atlanta}, \postcode{30332}, \state{Georgia}, \country{USA}}}

\affil[3]{\orgdiv{Department of Mathematics}, \orgname{University of Illinois Urbana-Champaign}, \orgaddress{\street{1409 W Green St}, \city{Urbana}, \postcode{61801}, \state{Illinois}, \country{USA}}}


\abstract{{\color{black} We study four synchronization notions for the all-to-all hybrid Kuramoto model containing both first- and second-order oscillators with heterogeneous inertias and damping coefficients.  We prove that full phase locking, bounded phase locking, and frequency synchronization are equivalent for arbitrary hybrid ensembles, and that each of these properties implies convergence of the complex order parameter.  Conversely, if an order-parameter synchronized trajectory has limiting order parameter $Z^*$ satisfying $|Z^*|>\max\left\{\frac{\omega_M}{\lambda},\,1-\frac{2}{N}\right\}$, then the trajectory is frequency synchronized and hence fully phase locked.  The converse proof is entirely real-dynamical.  The omega-limit set is internally chain transitive and, because the order parameter is constant on that set, the dynamics reduce there to a product of frozen scalar equations.  A tilted-energy argument shows that every nonstationary scalar factor of an internally chain transitive frozen set must cover the whole phase circle.  Evaluating the constant mean field at an antipodal phase then forces $|Z^*|\le 1-2/N$, a contradiction.  For zero natural frequencies, convergence follows from the analytic periodic gradient structure with degenerate inertia.}}

\keywords{First-, second-, and hybrid-order Kuramoto model; Synchronization; Order parameter}


\pacs[MSC Classification 2020]{Primary 34D45; Secondary 34D06, 34C25, 37C10}

\maketitle
\tableofcontents

\section{Introduction}
The Kuramoto model, one of the most representative and widely applied models of synchronization, was proposed nearly five decades ago (Kuramoto \cite{kuramoto1975self,kuramoto1984chemical}). Its synchronization-related applications span physics \cite{wiesenfeld1998frequency,yeung1999time, thumler2023synchrony,turner1998five,hsiao2023synchronization,hsiao2025synchronization}, biology \cite{bargiello1984restoration, hardin1990feedback, mirollo1990synchronization}, control theory \cite{bergen2007structure, filatrella2008analysis, dorfler2013synchronization,dorfler2012synchronization, li2014synchronization}, and non-Abelian and operator-valued extensions \cite{lohe2009non, bronski2020matrix,  hsiao2023synchronization, deville2019synchronization}, and we refer interested readers to \cite{strogatz2012sync} for a comprehensive overview. The equivalence of several synchronization notions for the first-order all-to-all Kuramoto flow was recently clarified in our related work~\cite{hsiaoLoZhu2026equivalence}. The corresponding question for genuinely hybrid systems, where first- and second-order oscillators coexist and inertia can support more complicated recurrent dynamics, is not covered by that result. The objective of the present work is to determine, under the simplest graph topology, namely all-to-all (complete) coupling with homogeneous interactions, which synchronization notions remain equivalent for the hybrid Kuramoto model with oscillator-dependent inertia and damping.

The study of synchronization has a rich and extensive history. Before briefly reviewing its development, let us consider the finite $N$ hybrid Kuramoto model consisting of both first- and second-order oscillators. The dynamics are given by  
\begin{equation}\label{mixed} 
\begin{aligned} 
d_j \dot{\theta}_j &= \omega_j + \frac{\lambda}{N} \sum_{k=1}^{N} \sin(\theta_k - \theta_j), \qquad j \in \{1,\ldots,n\}, \\[6pt] 
m_j \ddot{\theta}_j + d_j \dot{\theta}_j &= \omega_j + \frac{\lambda}{N} \sum_{k=1}^{N} \sin(\theta_k - \theta_j), \qquad j \in \{n+1,\ldots,N\}, 
\end{aligned} 
\end{equation}
where the inertia $m_j > 0$, damping $d_j > 0$, natural frequency $\omega_j \in \mathbb{R}$, and coupling strength $\lambda > 0$ are given. 
{\color{black} The parameter $n \in \{0,\ldots,N\}$ determines whether the system~\eqref{mixed} is purely first order ($n = N$), purely second order ($n = 0$), or genuinely hybrid.} 
The variable $\theta_j = \theta_j(t) \in \mathbb{R}$ denotes the phase of the $j$-th oscillator. Since the introduction of the order parameter by Kuramoto, measuring synchronization through this quantity has had a profound impact on scientific applications (see also Strogatz~\cite{strogatz2000kuramoto} and Acebrón et al.~\cite{acebron2005kuramoto}).
In this work, we refer to this notion as the \emph{order parameter synchronization state} (OPSS).
In addition to (OPSS), several alternative mathematical definitions are commonly employed to characterize and analyze synchronization, including the \emph{full phase-locked state} (FPLS), the \emph{phase-locked state} (PLS), and the \emph{frequency synchronization state} (FSS).

Previous studies have primarily focused on determining sufficient or necessary conditions for the emergence of synchronization under various parameter regimes of $m_j, d_j, \omega_j,$ and $\lambda$. 
{\color{black} The first-order Kuramoto model introduced in~\cite{kuramoto1975self,kuramoto1984chemical} is recovered from~\eqref{mixed} by taking $n = N$ and $d_j = 1$.} 
We refer the reader to~\cite{bronski2012fully, chopra2009exponential, dorfler2011critical, ha2020asymptotic, ha2016emergence, strogatz2000kuramoto, chen2024complete} for further developments on this classical setting. \textcolor{black}{More recently, Dai, Fiedler, and
López-Nieto~\cite{daiFiedlerLopezNieto2025} obtained a detailed global
description of the identical-frequency first-order Kuramoto flow by
establishing Morse--Smale structural stability and explicitly classifying
its heteroclinic connection graph.} For second-order Kuramoto models, the presence of inertia introduces substantially richer and more complex dynamics. 
Existing works have mainly established necessary or sufficient conditions for phase-locking or frequency synchronization in the inertial regime~\cite{ hsia2019synchronization, choi2011complete, choi2014complete, dorfler2011critical, cho2025inertia}. 
While these results provide criteria ensuring convergence to synchronized states, they typically address individual synchronization notions rather than the relation.

{\color{black} 
Clarifying the equivalence of these notions is important both mathematically and for applications. Phase-locking and frequency synchronization describe microscopic oscillator dynamics, whereas the order parameter provides a collective quantity that is often more convenient in physical applications and numerical studies. An equivalence theorem therefore allows synchronization detected through the order parameter to be rigorously interpreted in terms of the underlying phase and frequency dynamics, while also permitting one to use whichever characterization is most suitable for the analysis at hand.}

{\color{black} 
In this work we separate two issues that are often conflated.  First, for arbitrary
hybrid ensembles with oscillator-dependent inertias and damping coefficients, we prove
the unconditional equivalence
\[
    \mathrm{(FPLS)}
    \Longleftrightarrow
    \mathrm{(PLS)}
    \Longleftrightarrow
    \mathrm{(FSS)}.
\]
We denote these equivalent microscopic synchronization notions collectively by
$\mathrm{(SS)}$.  Every synchronized trajectory is also an order-parameter
synchronization state $\mathrm{(OPSS)}$ in the sense of Definition~\ref{def 3}.
The converse implication from the order parameter is more delicate because a frozen
second-order oscillator may possess nonstationary recurrent dynamics even when
equilibria exist.

Our second result gives a rigorous high-coherence converse.  If
$Z(t)\to Z^*$ and
\[
    |Z^*|>
    \max\left\{
       \frac{\omega_M}{\lambda},
       1-\frac{2}{N}
    \right\},
\]
then $\mathrm{(OPSS)}$ implies $\mathrm{(SS)}$.  The proof is purely
real-dynamical.  The omega-limit set of the original autonomous hybrid flow is compact
and internally chain transitive.  Since $Z$ is constant on this omega-limit set, the
restricted dynamics coincide exactly with the product of the frozen scalar equations
generated by the limiting mean field.  We show that a nonstationary scalar projection
of such an internally chain transitive set must cover the entire phase circle.  Indeed,
if one phase value were omitted, the scalar dynamics could be lifted to a compact
subset of the universal cover, where the tilted mechanical energy is a strict
Lyapunov function; the chain-transitivity convergence criterion of
Hirsch--Smith--Zhao~\cite{hirschSmithZhao2001}, together with the equivalence of
continuous-time and time-one chain recurrence~\cite{hurley1995chain}, then forces the
projected set to be an equilibrium.  Hence a genuinely nonstationary factor attains a
phase opposite to $Z^*$.  The constant-mean-field identity then gives the elementary
geometric obstruction
\[
    N|Z^*|\le N-2.
\]
This proves the stated threshold.

When all natural frequencies vanish, the stationary set need not be discrete.
In that case the hybrid system is an analytic periodic gradient-type system with
positive semidefinite inertia.  The convergence theorem of
Li--Xue~\cite[Theorem~1.1 and Remark~2.2]{liXue2019gradient} applies directly and
gives convergence of every trajectory to a stationary configuration.
Thus no finiteness assumption on the zero-frequency equilibrium set is required.}

\section{Definitions and Main Results} \label{sec 2}

We begin by introducing a co-rotating frame associated with a weighted average frequency. To avoid an unnecessary proliferation of notation, we shall, whenever convenient, write expressions such as $\omega_j\mapsto \omega_j-d_j(\omega_1+\ldots+\omega_N)/(d_1+\ldots+d_N)$, which is to be understood as ``$\omega_j$ is replaced by $\omega_j-d_j(\omega_1+\ldots+\omega_N)/(d_1+\ldots+d_N)$". With this convention, we perform the following change of variables:
\begin{align*}
\omega_j \mapsto \omega_j - \frac{\omega_1 + \cdots + \omega_N}{d_1 + \cdots + d_N} d_j, \quad \theta_j \mapsto \theta_j - \frac{\omega_1 + \cdots + \omega_N}{d_1 + \cdots + d_N} t,
\end{align*}
which in turn implies 
\begin{align*}
    \dot{\theta}_j \mapsto \dot{\theta}_j-\frac{\omega_1 + \cdots + \omega_N}{d_1 + \cdots + d_N}, \quad \ddot{\theta}_j\mapsto \ddot{\theta}_j\qquad\text{for all}\quad j\in \{1, \ldots, N\}.
\end{align*}
Under this change of frame, we may, without loss of generality, impose the normalization
\[
\omega_1 + \cdots + \omega_N = 0.
\]
which eliminates the uniform drift induced by the collective rotation. We remark that the special case $\omega_1=\omega_2=\ldots=\omega_N$ corresponds to identical natural frequencies. With this normalization in place, we now recall the standard definitions of synchronization in the Kuramoto model.

{\color{black} Throughout the remainder of the paper, we write
$\theta(t)=(\theta_1(t),\ldots,\theta_N(t))$ for the phase vector,
$\omega=(\omega_1,\ldots,\omega_N)$ for the natural-frequency vector, and
$\omega_M:=\max_{1\le j\le N}|\omega_j|$.}

\subsection{Definitions of Synchronization State}

\begin{definition}[PSS] \label{def 1}
A solution $\theta(t)$ is called a phase synchronization state if for all $j,k\in\{1,\ldots,N\}$,
    \begin{align*}
        \lim_{t\rightarrow \infty} \left(\theta_j(t)-\theta_k(t)\right)=0. 
    \end{align*}
\end{definition}

\begin{definition}[FPLS] \label{FPLS}
    A solution $\theta(t)$ is called a full phase-locked state if for all $j,k\in\{1,\ldots,N\}$,
    \begin{align*}
        \lim_{t\rightarrow\infty} (\theta_j(t)-\theta_k(t))=\theta^*_{jk},
    \end{align*}
where $\theta^*_{jk}$ are constants for all $j,k\in\{1,\ldots,N\}$.
\end{definition}

{\color{black}
\begin{remark}
Definition~\ref{FPLS} requires only convergence of all phase differences; it
does not assume a priori that the full phase vector converges.  Theorem~\ref{main 1 second}
proves the stronger conclusion that, in the co-rotating frame fixed above, an
(FPLS) trajectory converges to a stationary representative
$\theta^*=(\theta_1^*,\ldots,\theta_N^*)$.  In that case
$\theta^*_{jk}=\theta_j^*-\theta_k^*$, and $\theta^*$ satisfies
\[
    \omega_j+\frac{\lambda}{N}\sum_{k=1}^N
    \sin(\theta_k^*-\theta_j^*)=0,
    \qquad j=1,\ldots,N.
\]
\end{remark}
}

\begin{definition}[PLS]
    A solution $\theta(t)$ is called a phase-locked state if for all $j,k\in\{1,\ldots,N\}$, there exists a positive constant $L$ such that for all $t>0$,
    \begin{align*}
       \left|\theta_j(t)-\theta_k(t)\right|\leq L
       .
    \end{align*}
\end{definition}

\begin{definition}[FSS]
    A solution $\theta(t)$ is called a frequency synchronization state if 
    \begin{align*}
        \lim_{t\rightarrow \infty} \left|\dot{\theta}_j(t)\right|=0, ~~\mbox{for all}~~ j\in\{1,\ldots,N\}
        .
    \end{align*}
\end{definition}

{\color{black}
\begin{definition}[OPSS]\label{def 3}
For a solution $\theta(t)$, define the complex order parameter by
\begin{equation}\label{OP Z}
    Z(t):=\frac{1}{N}\sum_{j=1}^N e^{i\theta_j(t)}.
\end{equation}
We say that $\theta(t)$ is an order parameter synchronization state
$\mathrm{(OPSS)}$ if there exists $Z^*\in\mathbb C$ such that
\[
    Z(t)\longrightarrow Z^*
    \qquad\text{as }t\to\infty.
\]
\end{definition}

Thus (OPSS) refers only to convergence of the complex order parameter.  No
lower bound on $|Z^*|$ is included in the definition.  The estimate
$|Z^*|\ge\omega_M/\lambda$ will instead be derived below as a necessary
consequence of microscopic synchronization.  The converse implication requires
the additional strict high-coherence condition
\[
 |Z^*|>
 \max\left\{\frac{\omega_M}{\lambda},\,1-\frac{2}{N}\right\}.
\]
No claim is made here for the remaining low-coherence or boundary regime.
}
We pause to remark that a necessary condition for the existence of a phase synchronization state (PSS) is that all natural frequencies are identical. In other words, one must assume that  
\begin{align*}
    \omega_j -\omega_k \equiv 0, \quad \text{for all } j,k \in \{1, \dots, N\}.
\end{align*}  
Moreover, it is evident that if $\theta(t)$ is a full phase-locked state, then it is also a phase-locked state. {\color{black} Since the complex order parameter $Z(t)$ is defined directly by \eqref{OP Z}, no global choice of a continuous argument $\Theta(t)$ is required.} 

\subsection{Main Theorems} \label{sec 3}
The main result of this paper is on the equivalence of synchronization states for the hybrid Kuramoto model:
\begin{equation} \label{second main eq}
    \begin{aligned}
        &m_j\ddot{\theta}_j + d_j \dot{\theta}_j = \omega_j + \frac{\lambda}{N} \sum_{k=1}^{N} \sin(\theta_k - \theta_j), \qquad \sum_{j=1}^N \omega_j=0, \\
        &\qquad d_j>0,\quad \omega_j\in\mathbb{R},\quad \lambda>0, \qquad j \in \{1,\ldots,N\},\\
        &\qquad m_j=0, \quad j\in\{1,\ldots,n\}\quad\text{and}\quad m_j>0, \quad j\in\{n+1,\ldots,N\}.
    \end{aligned}
\end{equation}
{\color{black} We use the notation introduced above throughout the statement and proof of the main theorem.}
The main result of this paper is as follows.

{\color{black}
\begin{theorem}[Synchronization relations and a high-coherence converse]
\label{main 1 second}
Let $\theta(t)$ be a solution of the hybrid Kuramoto model
\eqref{second main eq}, with arbitrary $d_j>0$ and with
$m_j=0$ for the first-order components and $m_j>0$ for the second-order
components.
Then the following assertions hold.
\begin{enumerate}[label=\textnormal{(\roman*)}]
\item
The three microscopic synchronization notions are equivalent:
\[
    \mathrm{(FPLS)}
    \Longleftrightarrow
    \mathrm{(PLS)}
    \Longleftrightarrow
    \mathrm{(FSS)}.
\]
We denote any of these equivalent properties by $\mathrm{(SS)}$.

\item
Every synchronized trajectory is order-parameter synchronized:
\[
    \mathrm{(SS)}\Longrightarrow\mathrm{(OPSS)}.
\]
{\color{black} Moreover, if $Z(t)\to Z^*$ along a synchronized trajectory, then
\[
    |Z^*|\ge \frac{\omega_M}{\lambda}.
\]
Thus this lower bound is a necessary consequence of synchronization, rather
than part of the definition of (OPSS).}

\item
Conversely, assume that $\theta(t)$ is an $\mathrm{(OPSS)}$, {\color{black} so that}
\[
    Z(t)\longrightarrow Z^*
\]
and that either $\omega=0$, or
\begin{equation}\label{high coherence condition}
    |Z^*|
    >
    \max\left\{
       \frac{\omega_M}{\lambda},
       1-\frac{2}{N}
    \right\}.
\end{equation}
Then
\[
    \mathrm{(OPSS)}\Longrightarrow\mathrm{(SS)}.
\]
\end{enumerate}
In particular, in the nonzero-frequency case the four notions coincide along
every trajectory satisfying the strict high-coherence condition
\eqref{high coherence condition}.  No common-inertia or common-damping
assumption is imposed.
\end{theorem}
}

{\color{black}
\begin{remark}[Low-dimensional form of the high-coherence condition]
\label{remark low dimensional threshold}
For $\omega\neq0$, condition~\eqref{high coherence condition} becomes
\[
\begin{array}{ll}
N=2: & |Z^*|>\omega_M/\lambda,\\[1mm]
N=3: & |Z^*|>\max\{\omega_M/\lambda,\,1/3\},\\[1mm]
N=4: & |Z^*|>\max\{\omega_M/\lambda,\,1/2\}.
\end{array}
\]
Thus the present theorem does not assert an unconditional
$\mathrm{(OPSS)}\Rightarrow\mathrm{(SS)}$ implication for $N=3$ or $N=4$.
The zero-frequency case is separate: when $\omega=0$,
Lemma~\ref{zero frequency convergence} gives synchronization without any
additional coherence condition.
\end{remark}
}

{\color{black} 
The proof is organized in three parts.  We first establish the unconditional
equivalence among (FPLS), (PLS), and (FSS).  We then prove
$(\mathrm{SS})\Rightarrow(\mathrm{OPSS})$ directly from the limiting stationary
equations.  Finally, under \eqref{high coherence condition}, we derive the converse
from an internally chain transitive omega-limit-set argument and a scalar
tilted-energy rigidity lemma.}

\begin{lemma} \label{energy argument}
Let $\theta(t)$ be a solution of the hybrid Kuramoto model \eqref{second main eq}. 
If $\theta(t)$ is a $\mathrm{(PLS)}$, then it is also a $\mathrm{(FSS)}$.
\end{lemma}

\begin{proof}
Multiplying $\dot{\theta}_j$ in \eqref{second main eq}, summing over $j=1,\ldots,N$ and integrating the term over $(0,t)$, we obtain
\begin{equation*} \label{energy function}
    \begin{aligned}
        \sum_{j=1}^{N}\frac{m_j}{2}\dot{\theta}^2_j(t)+\int_0^t\sum_{j=1}^{N}d_j\dot{\theta}^2_j(s) \mathrm{d}s=&\sum_{j=1}^{N}\frac{m_j}{2}\dot{\theta}^2_j(0)+\sum_{j=1}^N \omega_j\left(\theta_j(t)-\theta_j(0)\right)\\
        &+\int_0^t \frac{\lambda}{N} \sum_{j=1}^N\sum_{k=1}^N\sin(\theta_k(s)-\theta_j(s))\dot{\theta}_j(s) \mathrm{d}s.
    \end{aligned}
    \end{equation*}

Since $\theta(t)$ is a (PLS), all phase differences $\theta_j-\theta_k$ remain uniformly bounded. Hence \begin{align*}
    \sum_{j=1}^N \omega_j\left(\theta_j(t)-\theta_j(0)\right)=\sum_{j=2}^N \omega_j(\theta_j(t)-\theta_1(t))-\sum_{j=1}^N \omega_j\theta_j(0)
\end{align*}
is uniformly bounded in $t$. Also, one may observe
    \begin{equation*} \label{H}
        \begin{aligned}
            \int_0^t \sum_{j=1}^N\sum_{k=1}^N\sin(\theta_k(s)-\theta_j(s))\dot{\theta}_j(s) \mathrm{d}s =\left.\sum_{j<k} \cos(\theta_k(s)-\theta_j(s))\right|_0^t,
        \end{aligned}
    \end{equation*}
    which is also uniformly bounded in $t$. Therefore, 
    \begin{align*}
        \lim\limits_{t\rightarrow \infty}\int_0^t\sum_{j=1}^{N}d_j\dot{\theta}^2_j(s) \mathrm{d}s\quad\text{exists and bounded}.
    \end{align*} 
    
{\color{black} By Lemma~\ref{a priori}, all velocities are uniformly bounded.
For $j>n$, the equation directly yields a uniform bound for
$\ddot\theta_j$.  For $j\le n$, differentiating the first-order equation gives
\[
 d_j\ddot\theta_j
 =
 \frac{\lambda}{N}\sum_{k=1}^N
 \cos(\theta_k-\theta_j)
 \bigl(\dot\theta_k-\dot\theta_j\bigr),
\]
so $\ddot\theta_j$ is uniformly bounded as well.  Hence
$\sum_j d_j\dot\theta_j^2$ has bounded derivative and is uniformly continuous.}
Applying Barbalat’s lemma \cite[Theorem 1]{farkas2016variations}, we obtain
\[
\lim_{t\to\infty}\sum_{j=1}^{N}d_j\dot{\theta}_j^2(t)=0,
\]
hence $\dot{\theta}_j(t)\to0$ for all $j$. Therefore, $\theta(t)$ is a (FSS).
\end{proof}

{\color{black} We proceed to prove the unconditional equivalence among (FPLS), (PLS), and (FSS).  The relation with (OPSS) is treated separately in Section~\ref{sec OPSS real}, where the converse implication is proved only under the high-coherence condition \eqref{high coherence condition}.} To streamline the argument, we begin with the following key lemma.

{\color{black}
\begin{lemma}[Zero-frequency convergence]\label{zero frequency convergence}
Assume $\omega_j=0$ for every $j$.  Then every solution of
\eqref{second main eq} converges to a stationary phase configuration:
there exists $\theta^*\in\mathbb R^N$ such that
\[
    \theta(t)\longrightarrow\theta^*,
    \qquad
    \dot\theta_j(t)\longrightarrow0
    \quad (j=1,\ldots,N).
\]
Consequently every solution is an (FPLS), a (PLS), an (FSS), and an (OPSS).
\end{lemma}

\begin{proof}
For $\omega=0$, define
\[
    \mathcal V(\theta)
    :=
    -\frac{\lambda}{N}
    \sum_{1\le j<k\le N}
    \cos(\theta_k-\theta_j).
\]
Then \eqref{second main eq} can be written
\[
    M\ddot\theta+D\dot\theta+\nabla\mathcal V(\theta)=0,
\]
where
\[
    M=\operatorname{diag}(m_1,\ldots,m_N)\ge0,
    \qquad
    D=\operatorname{diag}(d_1,\ldots,d_N)>0.
\]
The potential $\mathcal V$ is real analytic and $2\pi$-periodic in every
coordinate.  Li and Xue~\cite[Theorem~1.1]{liXue2019gradient} prove convergence
of every trajectory for analytic periodic first- and second-order gradient-type
systems; their Remark~2.2 states explicitly that the second-order result extends
to degenerate inertia
\[
    M=\operatorname{diag}(m_1,\ldots,m_q,0,\ldots,0)\ge0,
\]
which is precisely the mixed first-/second-order situation considered here.
Therefore every trajectory converges to a critical point of $\mathcal V$ and
$\dot\theta(t)\to0$.  Since the convergence asserted in that theorem is in
$\mathbb R^N$, it gives (FPLS) directly.  {\color{black} Finally, $Z(t)$ converges
automatically, so the trajectory is also an (OPSS).}
\end{proof}
}

{\color{black}
\begin{lemma}\label{finite root mode 2pi}
Assume $\omega_M>0$.  Let
$\psi=(\psi_1,\ldots,\psi_N)\in\mathbb T^N$ and consider
\begin{equation}\label{gj}
    g_j(\psi)
    :=
    \omega_j+\frac{\lambda}{N}
    \sum_{k=1}^{N}\sin(\psi_k-\psi_j)=0,
    \qquad j=1,\ldots,N,
\end{equation}
with $\sum_{j=1}^N\omega_j=0$.  Then the system either has no solution or
has only finitely many solutions modulo the common phase shift
$\psi\mapsto\psi+s\mathbf 1$.
\end{lemma}

\begin{proof}
Write
\[
    re^{i\Psi}:=\frac1N\sum_{k=1}^N e^{i\psi_k}.
\]
At a stationary point,
\[
    \omega_j+\lambda r\sin(\Psi-\psi_j)=0.
\]
Since $\omega_M>0$, every stationary point has
$r\ge\omega_M/\lambda>0$.  If $\lambda<\omega_M$, there is no stationary
point, so assume $\lambda\ge\omega_M$.  For each $j$,
\[
    \sin(\Psi-\psi_j)
    =
    -\frac{\omega_j}{\lambda r},
\]
and hence
\[
    \cos(\Psi-\psi_j)
    =
    \sigma_j
    \sqrt{1-\left(\frac{\omega_j}{\lambda r}\right)^2},
    \qquad \sigma_j\in\{-1,1\}.
\]
The self-consistency relation for the real part of the order parameter gives
\[
    rN
    =
    \sum_{j=1}^N
    \sigma_j
    \sqrt{1-\left(\frac{\omega_j}{\lambda r}\right)^2}.
\]
Set $x=\lambda r\in[\omega_M,\lambda]$.  For a fixed sign vector
$\sigma\in\{-1,1\}^N$, define
\[
    H_\sigma(x)
    :=
    x-\frac{\lambda}{N}
    \sum_{j=1}^N
    \sigma_j
    \sqrt{1-\left(\frac{\omega_j}{x}\right)^2},
    \qquad x\ge\omega_M.
\]
The function $H_\sigma$ is real analytic on $(\omega_M,\infty)$ and is not
identically zero there; indeed
\[
    H_\sigma(x)
    =
    x-\frac{\lambda}{N}\sum_{j=1}^N\sigma_j+O(x^{-2})
    \qquad (x\to\infty).
\]
It remains to exclude an accumulation of zeros at the endpoint
$x=\omega_M$.  Put $x=\omega_M+y^2$.  If
$|\omega_j|<\omega_M$, then
\[
    y\longmapsto
    \sqrt{1-
      \left(\frac{\omega_j}{\omega_M+y^2}\right)^2}
\]
is real analytic near $y=0$.  If $|\omega_j|=\omega_M$, then
\[
    \sqrt{1-
      \left(\frac{\omega_M}{\omega_M+y^2}\right)^2}
    =
    \frac{y\sqrt{2\omega_M+y^2}}{\omega_M+y^2},
    \qquad y\ge0,
\]
which is the restriction of a real-analytic function of $y$ near $0$.
Therefore
\[
    \widetilde H_\sigma(y)
    :=
    H_\sigma(\omega_M+y^2)
\]
extends real analytically across $y=0$.  If zeros of $H_\sigma$ accumulated
at $x=\omega_M$, the identity theorem for real-analytic functions would force
$\widetilde H_\sigma$, and hence $H_\sigma$, to vanish identically, a
contradiction.  Thus $H_\sigma$ has only finitely many zeros in
$[\omega_M,\lambda]$.

There are only finitely many sign vectors.  For each admissible pair
$(x,\sigma)$, the quantities
\[
    \sin(\Psi-\psi_j)
    =
    -\frac{\omega_j}{x},
    \qquad
    \cos(\Psi-\psi_j)
    =
    \sigma_j\sqrt{1-\left(\frac{\omega_j}{x}\right)^2}
\]
determine every relative angle $\Psi-\psi_j$ uniquely modulo $2\pi$.
Hence only finitely many relative-phase configurations occur, proving the
claim.
\end{proof}
}

We pause to remark that to the best of the authors’ knowledge, the finiteness of roots for such a system was first conjectured in~\cite{baillieul1982geometric}, where the case $N=3$ was rigorously proved and an intuitive argument for the finiteness when $N>3$ was also provided. For further studies on the finiteness of equilibria, we refer the reader to the numerical and formally analytical results in~\cite{mehta2015algebraic},~\cite{coss2018locating}, and~\cite{xi2017synchronization}, which confirm the finiteness of equilibria for various classes of coupling configurations and network topologies.

\begin{lemma}[Landau-Kolmogorov inequality]\label{Landau Kolmogorov}
{\color{black} Let $I=(a,\infty)$ with $a\ge0$, and let $f\in C^2(I)$ satisfy
$\|f\|_{L^\infty(I)}<\infty$ and $\|\ddot f\|_{L^\infty(I)}<\infty$. Then}
\[
{\color{black} \|\dot f\|_{L^\infty(I)}
\le 2\,\|f\|_{L^\infty(I)}^{1/2}
       \|\ddot f\|_{L^\infty(I)}^{1/2}.}
\]
\end{lemma}

\begin{proof}
{\color{black} Set $M:=\|f\|_{L^\infty(I)}$ and
$K:=\|\ddot f\|_{L^\infty(I)}$. For $x\in I$ and $h>0$, Taylor's theorem gives
\[
f(x+h)=f(x)+h\dot f(x)+\frac{h^2}{2}\ddot f(\xi)
\]
for some $\xi\in(x,x+h)$. Hence
\[
|\dot f(x)|\le \frac{2M}{h}+\frac{h}{2}K.
\]
If $M=0$ or $K=0$, the desired estimate is immediate. Otherwise, taking
$h=2\sqrt{M/K}$ yields
\[
|\dot f(x)|\le 2\sqrt{MK}.
\]
Taking the supremum over $x\in I$ proves the claim.}
\end{proof}
{\color{black} For $I=\mathbb{R}^+$ the constant $2$ is sharp
\cite{landau1925ungleichungen}; see also \cite{kolmogorov1949inequalities}
for more general Landau--Kolmogorov inequalities.}

Finally, we provide a fundamental a priori estimate for $\dot{\theta}(t)$.

\begin{lemma} \label{a priori}
The solution $\theta(t)$ of \eqref{second main eq} satisfies
\begin{equation*}
\begin{aligned}
\left|\dot{\theta}_j(t)\right| & \le \frac{\left|\omega_j\right|}{d_j}+\frac{\lambda}{d_j}, \qquad &j\in\{1,\ldots,n\},\\
\left|\dot{\theta}_j(t)\right| &\leq \left|\dot{\theta}_j(0)\right| e^{-\frac{d_j}{m_j} t} + {\color{black} \left(\frac{|\omega_j|}{d_j} + \frac{\lambda}{d_j}\right)}\left(1 - e^{-\frac{d_j}{m_j} t}\right), \qquad &j\in\{n+1,\ldots,N\},
\end{aligned}
\end{equation*}
for all $t \ge 0$.
\end{lemma}

{\color{black}
\begin{proof}
For $j\le n$, the first-order equation and
\[
 \left|\frac1N\sum_{k=1}^N
 \sin(\theta_k-\theta_j)\right|\le1
\]
give immediately
\[
    |\dot\theta_j(t)|
    \le\frac{|\omega_j|+\lambda}{d_j}.
\]
For $j>n$, set $v_j=\dot\theta_j$.  Variation of constants gives
\[
\begin{aligned}
v_j(t)
&=
e^{-d_jt/m_j}v_j(0)\\
&\quad
+\frac1{m_j}\int_0^t
e^{-d_j(t-s)/m_j}
\left[
 \omega_j+\frac{\lambda}{N}
 \sum_{k=1}^N\sin(\theta_k(s)-\theta_j(s))
\right]\,ds .
\end{aligned}
\]
Taking absolute values and evaluating the elementary exponential integral yields
\[
 |v_j(t)|
 \le
 |v_j(0)|e^{-d_jt/m_j}
 +
 \frac{|\omega_j|+\lambda}{d_j}
 \left(1-e^{-d_jt/m_j}\right),
\]
as claimed.
\end{proof}
}

\section{Proof of Theorem~\ref{main 1 second}}
We are now in a position to prove Theorem~\ref{main 1 second}.

\subsection{Proof of $\mathrm{(FPLS)}\Leftrightarrow \mathrm{(PLS)}\Leftrightarrow \mathrm{(FSS)}$}
\begin{lemma} \label{FSS implies PLS}
Let $\theta(t)$ be a solution of the hybrid Kuramoto model \eqref{second main eq}. 
Then the following implications hold:
\[
\mathrm{(FSS)} \Rightarrow \mathrm{(PLS)} 
\quad \text{and} \quad 
\mathrm{(PLS)} \Rightarrow \mathrm{(FPLS)}.
\]
\end{lemma}

\begin{proof}
(FSS) $\Rightarrow$ (PLS):
{\color{black} If $\omega_M=0$, then $\omega=0$ and
Lemma~\ref{zero frequency convergence} gives the stronger conclusion (FPLS).
Hence we may assume $\omega_M>0$.}
{\color{black} For every second-order component $j\in\{n+1,\ldots,N\}$,
the equation and Lemma~\ref{a priori} imply that $\ddot\theta_j$ is bounded;
differentiating the equation once then shows that $\dddot\theta_j$ is bounded as well.
Applying Lemma~\ref{Landau Kolmogorov} to
$f=\dot\theta_j$ on the half-line $(T,\infty)$ gives
\[
\|\ddot\theta_j\|_{L^\infty(T,\infty)}
\le
2\|\dot\theta_j\|_{L^\infty(T,\infty)}^{1/2}
\|\dddot\theta_j\|_{L^\infty(T,\infty)}^{1/2}.
\]
Since (FSS) gives
$\|\dot\theta_j\|_{L^\infty(T,\infty)}\to0$ as $T\to\infty$,
we conclude that $\ddot\theta_j(t)\to0$. Hence, if
\[
g(\theta):=(g_1(\theta),\ldots,g_N(\theta)),
\]
then the equations yield
\[
g_j(\theta(t))=
\begin{cases}
d_j\dot\theta_j(t), & j\le n,\\
m_j\ddot\theta_j(t)+d_j\dot\theta_j(t), & j>n,
\end{cases}
\]
and therefore
\[
\|g(\theta(t))\|\longrightarrow0.
\]}

{\color{black} We now make explicit the compactness argument.
Let
\[
\mathcal Q:=\mathbb T^N/\{\,\psi\sim\psi+s\mathbf 1\,\}
\cong\mathbb T^{N-1}
\]
be the relative-phase space, and let $q(t)=[\theta(t)]\in\mathcal Q$.
Because each $g_j$ depends only on phase differences, $g$ descends to a
continuous map $\bar g:\mathcal Q\to\mathbb R^N$.
The zero set of $\bar g$ is nonempty: by compactness of $\mathcal Q$, there
is a sequence $t_r\to\infty$ such that $q(t_r)\to q^\infty$, and the
convergence $g(\theta(t))\to0$ together with continuity gives
$\bar g(q^\infty)=0$.
By Lemma~\ref{finite root mode 2pi}, the zero set
\[
\mathcal E:=\{q\in\mathcal Q:\bar g(q)=0\}
=\{q^1,\ldots,q^M\}
\]
is finite. Choose $\epsilon>0$ sufficiently small so that the open balls
$B_\epsilon(q^\ell)\subset\mathcal Q$ are pairwise disjoint and are contained
in evenly covered neighborhoods for the universal covering
$\mathbb R^{N-1}\to\mathbb T^{N-1}\cong\mathcal Q$.
On the compact set
\[
K_\epsilon:=
\mathcal Q\setminus\bigcup_{\ell=1}^M B_\epsilon(q^\ell)
\]
there are no zeros of $\bar g$, and therefore
\[
\delta_\epsilon:=
\min_{q\in K_\epsilon}\|\bar g(q)\|>0.
\]
Since $\|g(\theta(t))\|\to0$, there exists $T>0$ such that
$\|g(\theta(t))\|<\delta_\epsilon$ for every $t\ge T$. Thus
\[
q(t)\in\bigcup_{\ell=1}^M B_\epsilon(q^\ell),
\qquad t\ge T.
\]
The set $q([T,\infty))$ is connected, because it is the continuous image of
the connected interval $[T,\infty)$. Since the above balls are pairwise
disjoint, the whole tail $q([T,\infty))$ must be contained in one fixed ball,
say $B_\epsilon(q^{\ell_0})$.}

{\color{black} Under the identification
$\mathcal Q\cong\mathbb T^{N-1}$, the corresponding lift is the relative-phase
vector
\[
\eta(t):=
(\theta_2(t)-\theta_1(t),\ldots,\theta_N(t)-\theta_1(t))
\in\mathbb R^{N-1}.
\]
Because $B_\epsilon(q^{\ell_0})$ is evenly covered and $\eta(t)$ is continuous,
the lifted tail $\eta([T,\infty))$ remains in one connected component of the
preimage of this ball. That component is bounded. Consequently all relative
phases $\theta_j(t)-\theta_k(t)$ are uniformly bounded, and hence
$\theta(t)$ is a (PLS).}

(PLS) $\Rightarrow$ (FPLS):
{\color{black} If $\omega_M=0$, the conclusion follows directly from
Lemma~\ref{zero frequency convergence}.  Hence assume $\omega_M>0$.}
{\color{black} By Lemma~\ref{energy argument}, a (PLS) is also a (FSS), and
the preceding Landau--Kolmogorov argument again gives
$\ddot\theta_j(t)\to0$ for all second-order components. Summing
\eqref{second main eq} over $j$ and using
$\sum_j\omega_j=0$ together with the cancellation of the coupling terms yields
the exact identity
\[
\sum_{j=1}^N m_j\dot\theta_j(t)
+\sum_{j=1}^N d_j\theta_j(t)
=
\sum_{j=1}^N m_j\dot\theta_j(0)
+\sum_{j=1}^N d_j\theta_j(0).
\]
Since $\dot\theta_j(t)\to0$, we obtain
\begin{align}\label{g0}
g_0(\theta(t))
:=\sum_{j=1}^N d_j\theta_j(t)
\longrightarrow
C_0
:=
\sum_{j=1}^N m_j\dot\theta_j(0)
+\sum_{j=1}^N d_j\theta_j(0).
\end{align}
This is the role of the constraint $g_0$: it controls the neutral common
phase-shift direction. Indeed,
\[
g_0(\theta(t))
=
\Big(\sum_{j=1}^N d_j\Big)\theta_1(t)
+\sum_{j=2}^N d_j\big(\theta_j(t)-\theta_1(t)\big).
\]
The second term is bounded by (PLS), while $g_0(\theta(t))$ converges.
Therefore $\theta_1(t)$ is bounded, and hence all $\theta_j(t)$ are bounded.}

{\color{black} Let $\Omega$ denote the $\omega$-limit set of the bounded
trajectory $\theta(t)$. It is nonempty and compact. Moreover,
\[
\Omega=\bigcap_{T>0}\overline{\{\theta(t):t\ge T\}},
\]
and the sets on the right are nested compact connected sets; hence $\Omega$
is connected. Since $\dot\theta_j(t)\to0$ and
$\ddot\theta_j(t)\to0$ for the second-order components, every
$\theta^*\in\Omega$ satisfies
\[
g_j(\theta^*)=0,\qquad j=1,\ldots,N,
\qquad\text{and}\qquad
g_0(\theta^*)=C_0.
\]
By Lemma~\ref{finite root mode 2pi}, there are only finitely many stationary
configurations modulo a common phase shift. Intersecting these classes with
the condition $g_0=C_0$ removes the common-shift freedom; since the trajectory
is bounded, only finitely many lifted representatives can occur. Thus
$\Omega$ is contained in a finite set. A connected subset of a finite set is
a singleton, say $\Omega=\{\theta^*\}$. Therefore
\[
\theta(t)\longrightarrow\theta^*,
\]
so $\theta(t)$ is a (FPLS).}
This completes the proof of Lemma~\ref{FSS implies PLS}.
\end{proof}

Combining Lemma~\ref{energy argument} with Lemma~\ref{FSS implies PLS}, we obtain the equivalence among (FPLS), (PLS), and (FSS). From now on, we refer to these equivalent notions collectively as the synchronization state (SS). 

\subsection{{\color{black} Order-parameter synchronization and the high-coherence converse}}
\label{sec OPSS real}

\begin{lemma}\label{SS implies OPSS}
Let $\theta(t)$ be a solution of \eqref{second main eq}.  Then
\[
    \mathrm{(SS)}\Longrightarrow\mathrm{(OPSS)}.
\]
\end{lemma}

{\color{black}
\begin{proof}
By the equivalence already proved, an (SS) is both an (FPLS) and an (FSS).
Hence there exist $\theta_j^*\in\mathbb R$ such that
\[
    \theta_j(t)\to\theta_j^*,
    \qquad
    \dot\theta_j(t)\to0.
\]
For every second-order component, the Landau--Kolmogorov argument used in
Lemma~\ref{FSS implies PLS} also gives $\ddot\theta_j(t)\to0$.  Therefore
\[
    Z(t)\longrightarrow
    Z^*:=
    \frac1N\sum_{j=1}^N e^{i\theta_j^*}.
\]
Passing to the limit in the equations gives
\[
    \omega_j+
    \lambda\operatorname{Im}
    \bigl(Z^*e^{-i\theta_j^*}\bigr)=0,
    \qquad j=1,\ldots,N.
\]
Consequently
\[
    |\omega_j|\le\lambda|Z^*|
    \qquad\text{for every }j,
\]
and hence $|Z^*|\ge\omega_M/\lambda$.  {\color{black} Thus the solution is an
(OPSS), and the stated lower bound follows as an additional necessary
consequence of synchronization.}
\end{proof}
}

We next establish the real-dynamical rigidity statement used for the converse.

{\color{black}
\begin{lemma}[Scalar phase coverage for internally chain transitive frozen sets]
\label{scalar ICT phase coverage}
Let $a>|\omega|$, $d>0$, and let either $m=0$ or $m>0$.  Consider the
frozen scalar equation on the phase circle
\begin{equation}\label{scalar frozen real}
    m\ddot x+d\dot x=\omega-a\sin x,
\end{equation}
with the usual first-order interpretation when $m=0$.  Let $K$ be a nonempty
compact invariant internally chain transitive set for its continuous flow.
If $K$ is not a single equilibrium, then the phase projection of $K$ is the
whole circle:
\[
    \pi_x(K)=\mathbb T.
\]
\end{lemma}

\begin{proof}
Assume, to the contrary, that $\pi_x(K)\neq\mathbb T$.  Choose
$\beta\in\mathbb T\setminus\pi_x(K)$.  The quotient map restricts to a
homeomorphism from $(\beta,\beta+2\pi)$ onto $\mathbb T\setminus\{\beta\}$,
so its inverse gives a single continuous lift of the phase coordinate on
$K$.  Since $\pi_x(K)$ is compact and avoids $\beta$, there is a
$\delta>0$ such that the lifted phase satisfies
\[
    x\in[\beta+\delta,\,\beta+2\pi-\delta]
    \qquad\text{on }K.
\]
Lifting the velocity coordinate unchanged when $m>0$ gives a compact set
$\widetilde K$ homeomorphic to $K$.  Because an orbit in $K$ never reaches
the cut point $\beta$, its uniquely lifted orbit remains in the same branch;
hence $\widetilde K$ is invariant for the lifted scalar flow.  In particular,
the energy introduced below is a single-valued $C^1$ function on the compact
invariant set $\widetilde K$.

For $m>0$ define the lifted tilted mechanical energy
\[
    \mathcal E(x,v)
    :=
    \frac{m}{2}v^2-\omega x-a\cos x,
    \qquad v=\dot x.
\]
For $m=0$ set
\[
    \mathcal E(x):=-\omega x-a\cos x.
\]
Along the scalar flow one has in both cases
\begin{equation}\label{scalar energy dissipation}
    \frac{d}{dt}\mathcal E
    =
    -d\,\dot x^{\,2}\le0,
\end{equation}
with equality along an entire orbit only at an equilibrium.

Because $a>|\omega|$, the equilibria of \eqref{scalar frozen real} are
nondegenerate and therefore isolated.  On the compact invariant set
$\widetilde K$,
\eqref{scalar energy dissipation} and LaSalle's invariance principle imply
that every positive orbit has an omega-limit set consisting of equilibria;
since the equilibria are isolated, every orbit in $K$ converges to one
equilibrium.  Moreover, there is no cycle of equilibria: along every
nonconstant full connecting orbit the value of $\mathcal E$ decreases
strictly.

Every point of an internally chain transitive set is chain recurrent for the
restricted scalar semiflow.  Hurley~\cite{hurley1995chain} proves that the
chain recurrent set of a continuous semiflow agrees with that of its time-one
map.  Hence every $y\in K$ is chain recurrent for the time-one map
$f:=\Phi(1)|_K$.  Let $C_y$ be the chain component of $y$ for $f$.  Since
the chain relation of a continuous map on a compact metric space is closed,
$C_y$ is closed in $K$ and hence compact; mutual chain reachability within
the component also makes $C_y$ internally chain transitive for $f$.

We verify the hypotheses of the convergence theorem of Hirsch, Smith, and
Zhao~\cite[Theorem~3.2]{hirschSmithZhao2001} for $f$.  First, a fixed point
of $f$ must be an equilibrium of the continuous scalar flow: otherwise
\eqref{scalar energy dissipation} would give
$\mathcal E(\Phi(1)y)<\mathcal E(y)$, contradicting $f(y)=y$.
The strict inequality $a>|\omega|$ makes all scalar equilibria hyperbolic,
hence isolated invariant sets.  Second, every orbit in the compact set $K$
converges to an equilibrium.  Indeed, LaSalle's invariance principle applied
to \eqref{scalar energy dissipation} puts every omega-limit set in the finite
set of equilibria contained in the lifted compact strip; an omega-limit set
of a continuous orbit is connected, and therefore consists of one
equilibrium.  Third, no cycle of fixed points is possible, because along
every nonconstant connecting orbit the value of $\mathcal E$ decreases
strictly.

Theorem~3.2 of~\cite{hirschSmithZhao2001} therefore implies that the compact
internally chain transitive set $C_y$ is a single fixed point.  Thus every
$y\in K$ is an equilibrium of the scalar flow.  Since the equilibria are
isolated and $K$ itself is internally chain transitive, $K$ must be a single
equilibrium, contradicting the hypothesis.  Therefore
$\pi_x(K)=\mathbb T$.
\end{proof}

\begin{lemma}[High-coherence rigidity of the frozen mean field]
\label{high coherence frozen rigidity}
Let $R\in(0,1]$ and $a=\lambda R$, and suppose
\[
    a>\omega_M.
\]
Let $\Omega$ be a nonempty compact invariant internally chain transitive set
for the frozen hybrid product flow
\begin{equation}\label{real frozen product}
    m_j\ddot x_j+d_j\dot x_j
    =
    \omega_j-a\sin x_j,
    \qquad j=1,\ldots,N,
\end{equation}
and assume that
\begin{equation}\label{real frozen constant mean}
    \frac1N\sum_{j=1}^N e^{ix_j}=R
    \qquad\text{on }\Omega.
\end{equation}
If $\Omega$ contains a non-equilibrium point, then
\begin{equation}\label{coherence obstruction}
    R\le1-\frac{2}{N}.
\end{equation}
\end{lemma}

\begin{proof}
Assume that $\Omega$ contains a non-equilibrium point.  Then for some index
$j$ the scalar projection
\[
    K_j:=\pi_j(\Omega)
\]
is not a single equilibrium.  The projection of an internally chain
transitive set under a continuous factor map is internally chain transitive;
hence $K_j$ is a compact invariant internally chain transitive set for the
$j$-th scalar flow.  Since
\[
    a>\omega_M\ge|\omega_j|,
\]
Lemma~\ref{scalar ICT phase coverage} yields
\[
    \pi_{x_j}(K_j)=\mathbb T.
\]
Therefore there exists a point of $\Omega$ for which
\[
    e^{ix_j}=-1.
\]
Taking the real part of \eqref{real frozen constant mean} at this point gives
\[
    NR
    =
    -1+\sum_{k\ne j}\cos x_k
    \le
    -1+(N-1)
    =
    N-2.
\]
This is exactly \eqref{coherence obstruction}.
\end{proof}

\begin{lemma}[High-coherence OPSS implies synchronization]
\label{OPSS implies SS}
Let $\theta(t)$ be a solution of \eqref{second main eq} which is an (OPSS),
with
\[
    Z(t)\longrightarrow Z^*.
\]
If $\omega=0$, or if
\[
    |Z^*|
    >
    \max\left\{
       \frac{\omega_M}{\lambda},
       1-\frac{2}{N}
    \right\},
\]
then $\theta(t)$ is an (SS).
\end{lemma}

\begin{proof}
If $\omega=0$, the conclusion follows from
Lemma~\ref{zero frequency convergence}.  Assume therefore $\omega\neq0$ and
the strict high-coherence condition above.

The original hybrid Kuramoto system is autonomous on
\[
    \mathcal X
    :=
    \mathbb T^N\times\mathbb R^{N-n},
\]
where the second factor consists of the velocities of the inertial
components.  Lemma~\ref{a priori} makes the positive orbit precompact.
Let $\Omega$ be its omega-limit set.  By the continuous-semiflow version of
the omega-limit theorem, $\Omega$ is nonempty, compact, invariant, and
internally chain transitive; see
Hirsch--Smith--Zhao~\cite[Lemma~2.1$'$]{hirschSmithZhao2001}.

Because $Z(t)\to Z^*$ and $Z$ is continuous,
\[
    Z\equiv Z^*
    \qquad\text{on }\Omega.
\]
Before freezing, perform one common phase rotation in the original
all-to-all system so that
\[
    Z^*=R^*>0.
\]
This is a genuine symmetry of the original Kuramoto equations because the
coupling depends only on phase differences.  Since $\Omega$ is invariant,
every complete orbit contained in $\Omega$ has order parameter identically
equal to $R^*$.  Consequently, on $\Omega$ the original vector field agrees
exactly with the frozen product flow
\[
    m_j\ddot x_j+d_j\dot x_j
    =
    \omega_j-\lambda R^*\sin x_j.
\]
Thus $\Omega$ satisfies the hypotheses of
Lemma~\ref{high coherence frozen rigidity}, with
$a=\lambda R^*>\omega_M$.

If $\Omega$ contained a non-equilibrium point, that lemma would imply
\[
    R^*\le1-\frac2N,
\]
contrary to the assumed high-coherence condition.  Hence every point of
$\Omega$ is an equilibrium of the frozen product flow.

It follows that the original trajectory is an (FSS).  Indeed, if a
second-order velocity failed to converge to zero, a subsequence of states
would converge to a point of $\Omega$ with a nonzero velocity component,
contradicting the preceding conclusion.  For a first-order component,
\[
    d_j\dot\theta_j
    =
    \omega_j+
    \lambda\operatorname{Im}
    \bigl(Z(t)e^{-i\theta_j}\bigr)
\]
is a continuous function of the phase state and $Z(t)$; a subsequential
nonzero limiting velocity would again give a non-equilibrium point of
$\Omega$.  Therefore
\[
    \dot\theta_j(t)\longrightarrow0
    \qquad (j=1,\ldots,N).
\]
By Lemma~\ref{FSS implies PLS}, (FSS) is equivalent to (SS), completing the
proof.
\end{proof}
}

Finally, we are ready to prove Theorem~\ref{main 1 second}.

{\color{black}
\begin{proof}[Proof of Theorem~\ref{main 1 second}]
Lemma~\ref{energy argument} and Lemma~\ref{FSS implies PLS}, together with
the trivial implication $(\mathrm{FPLS})\Rightarrow(\mathrm{PLS})$, give
\[
    \mathrm{(FPLS)}
    \Longleftrightarrow
    \mathrm{(PLS)}
    \Longleftrightarrow
    \mathrm{(FSS)}.
\]
Lemma~\ref{SS implies OPSS} proves
$\mathrm{(SS)}\Rightarrow\mathrm{(OPSS)}$.
Finally, Lemma~\ref{OPSS implies SS} proves the converse in the
zero-frequency case and under the strict high-coherence condition
\eqref{high coherence condition}.  This proves all assertions.
\end{proof}
}

\color{black}
\section{Numerical simulations}
\label{sec:numerics}

In this section, we present numerical simulations illustrating the
synchronization relations established in Theorem~\ref{main 1 second} for genuinely hybrid
ensembles containing both first- and second-order oscillators. Our purpose is
not to determine a numerical critical coupling, but rather to examine whether
the microscopic and macroscopic synchronization diagnostics appearing in the
theorem exhibit the predicted concurrent asymptotic behavior under
heterogeneous inertia, damping, and time scales.

We consider the hybrid Kuramoto system
\begin{equation}
m_j \ddot{\theta}_j + d_j \dot{\theta}_j
=
\omega_j
+
\frac{\lambda}{N}
\sum_{k=1}^{N}
\sin(\theta_k-\theta_j),
\qquad
j=1,\ldots,N,
\label{eq:num-hybrid-system}
\end{equation}
with the first-order interpretation when $m_j=0$. As in the theoretical
analysis, the natural frequencies are centered so that
\begin{equation}
\sum_{j=1}^{N} \omega_j = 0.
\label{eq:num-centered-frequency}
\end{equation}

Because strongly heterogeneous inertia and damping may generate stiff
transient dynamics, we use MATLAB's \texttt{ode15s} solver. For the principal
trajectory simulations, the relative and absolute error tolerances are
$10^{-8}$ and $10^{-10}$, respectively. The fine coupling sweep described
below is recomputed using the tighter tolerances $10^{-9}$ and $10^{-11}$.
The complete realized parameter vectors, initial states, random seed, solver
options (including maximum step sizes), terminal integration intervals, and
post-processing scripts underlying the reported diagnostics are available as
stated in the Data and Code Availability Statement.

To quantify the terminal behavior, let $I_{\rm tail}$ denote the terminal
portion of a numerical trajectory. We record
\begin{equation}
V_{\rm tail}
:=
\max_{t\in I_{\rm tail}}
\max_{1\leq j\leq N}
\left|\dot{\theta}_j(t)\right|,
\label{eq:num-vtail}
\end{equation}
as a finite-time frequency-synchronization diagnostic. We also define
\begin{equation}
D_Z
:=
\max_{t\in I_{\rm tail}}
\left|
Z(t)-\overline{Z}_{\rm tail}
\right|,
\label{eq:num-dz}
\end{equation}
where $\overline{Z}_{\rm tail}$ denotes the numerical mean of $Z(t)$ over
$I_{\rm tail}$.

To monitor convergence of relative phases, we divide $I_{\rm tail}$ into its
first and second halves, denoted by $I_{\rm tail}^{(1)}$ and
$I_{\rm tail}^{(2)}$, and define
\begin{equation}
D_{\rm rel}
:=
\max_{2\leq j\leq N}
\left|
\left\langle \theta_j-\theta_1 \right\rangle_{I_{\rm tail}^{(2)}}
-
\left\langle \theta_j-\theta_1 \right\rangle_{I_{\rm tail}^{(1)}}
\right|,
\label{eq:num-drel}
\end{equation}
where $\langle\cdot\rangle_I$ denotes the numerical average over the time
interval $I$. These quantities are used only as finite-time numerical
diagnostics and are not substitutes for the asymptotic definitions introduced
in Section~2.

\subsection{Paired strong- and weak-coupling hybrid dynamics}
\label{subsec:num-strong-weak}

We first consider a genuinely hybrid ensemble with
\begin{equation}
N=30,
\qquad
n=15.
\label{eq:num-N30}
\end{equation}
A single realization of the natural frequencies, inertias, damping
coefficients, and initial conditions is generated and then held fixed. The
natural frequencies are sampled from a Gaussian distribution with standard
deviation $0.35$ and are subsequently centered according to
\eqref{eq:num-centered-frequency}. For the realized sample,
\begin{equation}
\omega_M \approx 0.8800.
\label{eq:num-N30-omegaM}
\end{equation}
For the second-order components, the positive inertias satisfy approximately
\begin{equation}
0.1027 \leq m_j \leq 9.8538,
\label{eq:num-N30-mrange}
\end{equation}
while the damping coefficients of the full ensemble satisfy approximately
\begin{equation}
0.3375 \leq d_j \leq 2.9583.
\label{eq:num-N30-drange}
\end{equation}

\begin{figure}[t]
    \centering
    \includegraphics[width=\textwidth]{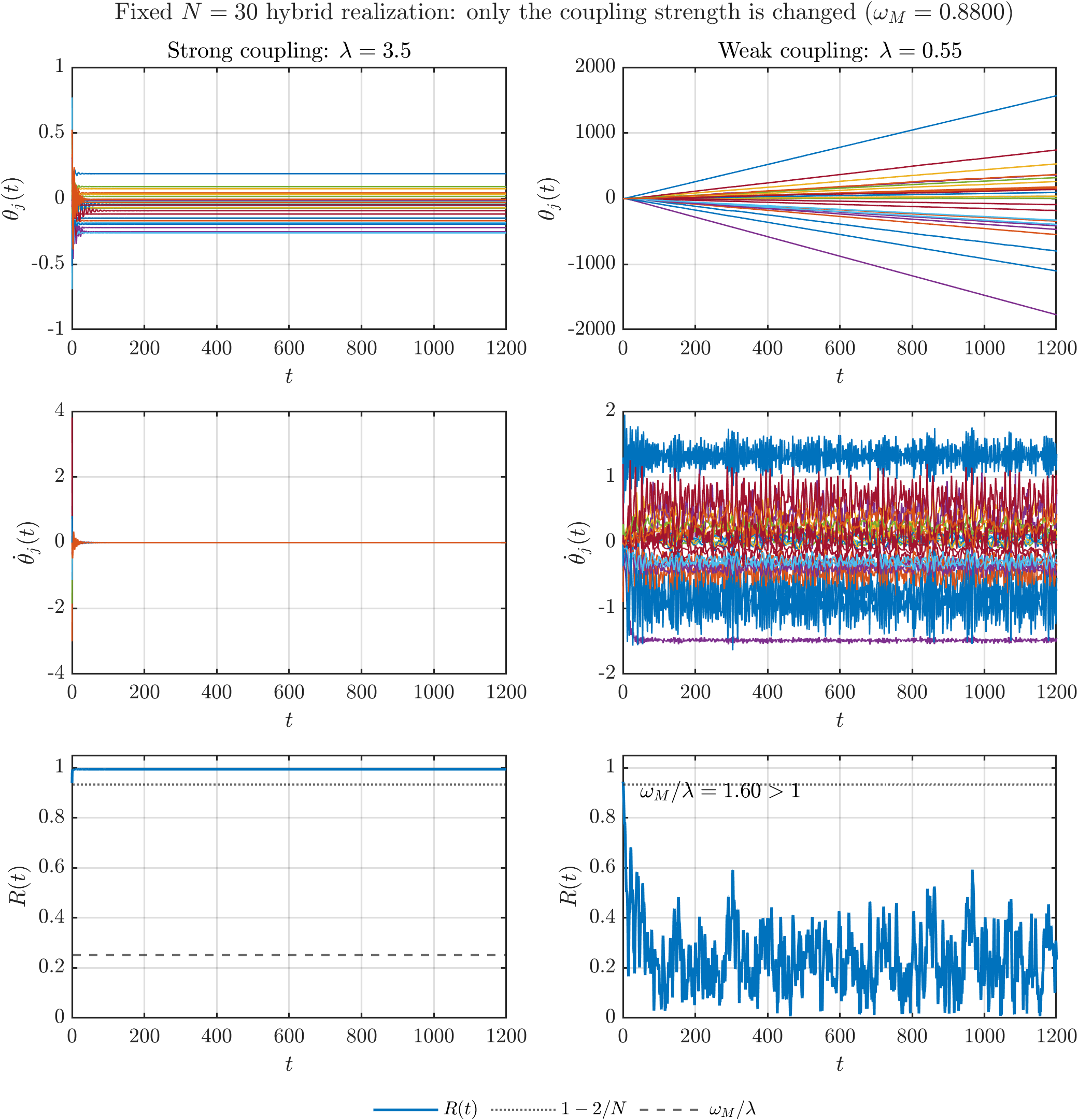}
    \caption{
    Paired simulations of the same genuinely hybrid $N=30$ realization with
    the same realized natural-frequency vector, inertia coefficients, damping
    coefficients, and initial state, but different coupling strengths. For
    $\lambda=3.5$ (left), the phases lock, all frequencies decay to zero, and
    the order parameter converges. For $\lambda=0.55$ (right), the phases
    continue to drift, the frequencies remain nonzero, and the order parameter
    does not converge. The comparison isolates the effect of the coupling
    strength and illustrates the concurrent behavior of the microscopic and
    macroscopic synchronization diagnostics.
    }
    \label{fig:hybrid-strong-weak}
\end{figure}

Figure~\ref{fig:hybrid-strong-weak} compares two trajectories for which the
only changed parameter is the coupling strength. For strong coupling,
$\lambda=3.5$, the trajectory converges numerically to a strongly coherent
synchronized configuration. The terminal order-parameter magnitude is
\begin{equation}
R_{\rm tail} \approx 0.994438,
\label{eq:num-strong-R}
\end{equation}
while the terminal diagnostics are
\begin{align}
V_{\rm tail}
&\approx 6.13\times 10^{-16},
\nonumber\\
D_{\rm rel}
&\approx 5.55\times 10^{-17},
\nonumber\\
D_Z
&\approx 8.95\times 10^{-16}.
\label{eq:num-strong-diagnostics}
\end{align}

In contrast, when the coupling is reduced to
\begin{equation}
\lambda=0.55,
\label{eq:num-weak-lambda}
\end{equation}
while every other parameter and the initial condition are left unchanged,
the trajectory remains nonsynchronized. Numerically, we obtain
\begin{align}
R_{\rm tail}
&\approx 0.227391,
\nonumber\\
\operatorname{std}_{I_{\rm tail}} R(t)
&\approx 0.110245,
\label{eq:num-weak-R}
\end{align}
together with
\begin{align}
V_{\rm tail}
&\approx 1.69,
\nonumber\\
D_{\rm rel}
&\approx 2.52\times 10^{2},
\nonumber\\
D_Z
&\approx 0.514.
\label{eq:num-weak-diagnostics}
\end{align}

Thus, for this fixed hybrid ensemble, the phase, frequency, and
order-parameter diagnostics either stabilize concurrently or fail to stabilize
concurrently.

\subsection{Fine coupling sweep and the high-coherence condition}
\label{subsec:num-lambda-sweep}

To examine more systematically the relation between synchronization and the
strict high-coherence condition, we keep the entire $N=30$ realization from
Section~\ref{subsec:num-strong-weak} fixed and vary only the coupling strength
over the grid
\begin{equation}
\lambda
\in
\{0.90,0.91,\ldots,1.20\}.
\label{eq:num-lambda-grid}
\end{equation}
Each trajectory is integrated to $T=2500$.

For this realization,
\begin{equation}
1-\frac{2}{N}
=
1-\frac{2}{30}
=
0.933333\ldots .
\label{eq:num-N30-geom-threshold}
\end{equation}

\begin{figure}[t]
    \centering
    \includegraphics[width=\textwidth]{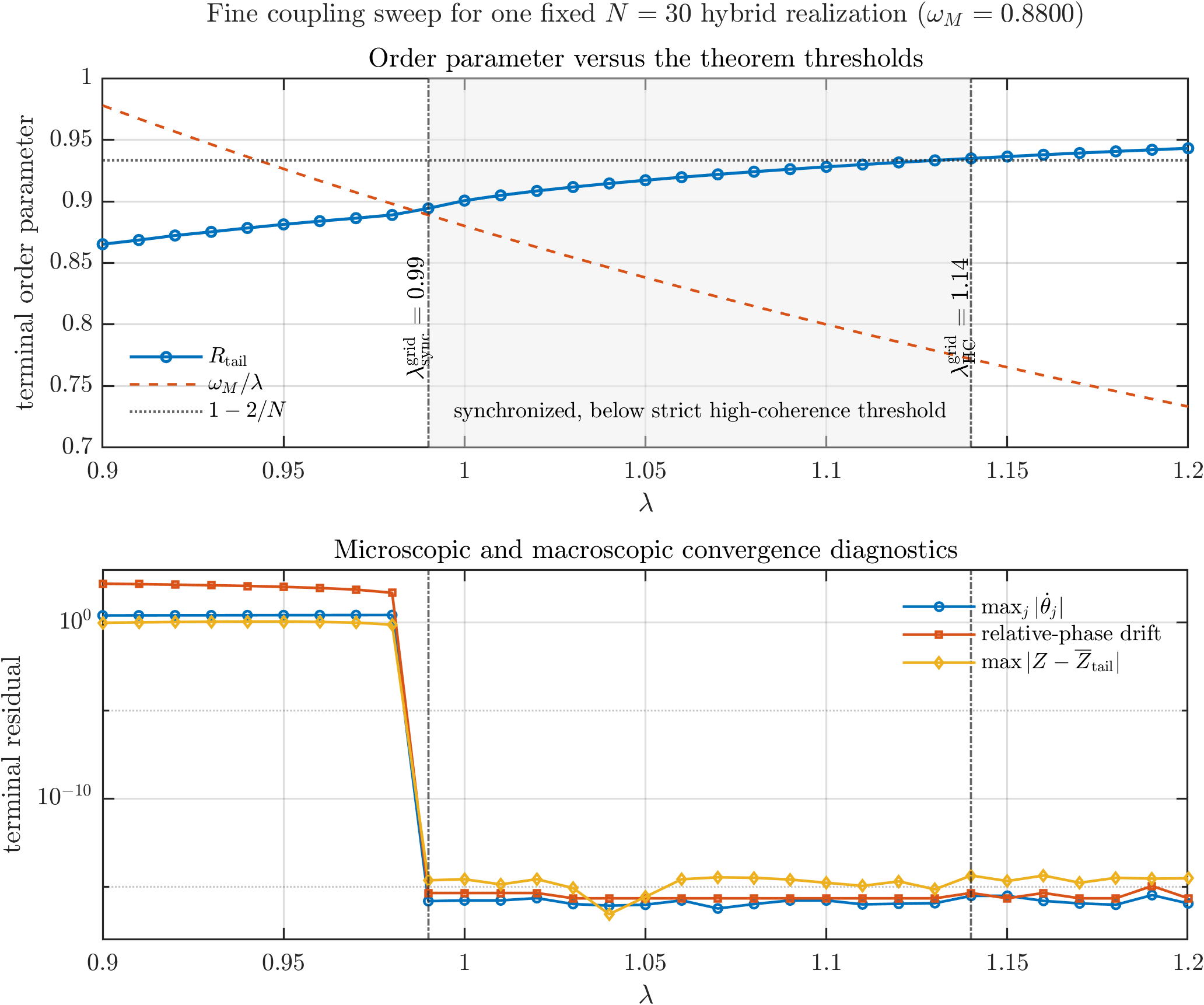}
    \caption{
    Fine coupling sweep for one fixed genuinely hybrid $N=30$ realization.
    The top panel compares the terminal order-parameter magnitude with
    $\omega_M/\lambda$ and $1-2/N$. The shaded interval denotes sampled
    coupling strengths for which the trajectory is already numerically
    synchronized but remains below the strict high-coherence threshold. The
    bottom panel shows the simultaneous collapse of the frequency,
    relative-phase, and complex-order-parameter residuals. The figure
    emphasizes that the high-coherence condition is a sufficient hypothesis
    for the converse theorem rather than a synchronization-onset criterion.
    }
    \label{fig:lambda-sweep}
\end{figure}

The resulting numerical diagnostics are shown in
Figure~\ref{fig:lambda-sweep}. On the sampled grid, the trajectory is clearly
nonsynchronized at $\lambda=0.98$. In particular,
\begin{align}
R_{\rm tail}
&\approx 0.888949,
\nonumber\\
V_{\rm tail}
&\approx 2.64,
\nonumber\\
D_{\rm rel}
&\approx 48.5,
\nonumber\\
D_Z
&\approx 0.753.
\label{eq:num-lambda098}
\end{align}

At the next sampled value, $\lambda=0.99$, all three convergence residuals
collapse to numerical precision:
\begin{align}
R_{\rm tail}
&\approx 0.894318,
\nonumber\\
V_{\rm tail}
&\approx 1.53\times 10^{-16},
\nonumber\\
D_{\rm rel}
&\approx 4.44\times 10^{-16},
\nonumber\\
D_Z
&\approx 2.34\times 10^{-15}.
\label{eq:num-lambda099}
\end{align}

Moreover, at $\lambda=0.99$,
\begin{equation}
\frac{\omega_M}{\lambda}
\approx 0.888925,
\label{eq:num-lambda099-freqthreshold}
\end{equation}
and hence
\begin{equation}
R_{\rm tail}
>
\frac{\omega_M}{\lambda}.
\label{eq:num-lambda099-opss}
\end{equation}
Thus, the numerical trajectory satisfies the frequency-amplitude inequality
$R_{\rm tail}>\omega_M/\lambda$ appearing in the sufficient high-coherence
hypothesis of Theorem~\ref{main 1 second}; this inequality is not part of the
definition of order-parameter synchronization.  However,
\begin{equation}
R_{\rm tail}
<
1-\frac{2}{N}
=
0.933333\ldots .
\label{eq:num-lambda099-not-hc}
\end{equation}

The geometric high-coherence threshold is not crossed until between
$\lambda=1.13$ and $\lambda=1.14$. Indeed,
\begin{align}
R_{\rm tail}(1.13)
&\approx 0.933271
<
0.933333\ldots,
\nonumber\\
R_{\rm tail}(1.14)
&\approx 0.934859
>
0.933333\ldots .
\label{eq:num-hc-crossing}
\end{align}

Accordingly, on the sampled grid we may summarize the two numerical
landmarks as
\begin{equation}
\lambda_{\rm sync}^{\rm grid}=0.99,
\qquad
\lambda_{\rm HC}^{\rm grid}=1.14.
\label{eq:num-grid-landmarks}
\end{equation}
Here $\lambda_{\rm sync}^{\rm grid}$ denotes the first sampled coupling for
which the numerical synchronization diagnostics fall below the prescribed
terminal tolerances, while $\lambda_{\rm HC}^{\rm grid}$ denotes the first
sampled coupling for which the strict high-coherence condition is satisfied.
Neither quantity is asserted to be an analytically defined critical coupling.

The interval between these two numerical landmarks provides a useful
illustration of the logical role of the high-coherence assumption in
Theorem~\ref{main 1 second}: it is a sufficient condition under which order-parameter
synchronization implies microscopic synchronization, but it is not asserted
to be a necessary condition for synchronization itself.

\subsection{Zero natural frequencies with extreme scale separation}
\label{subsec:num-zero-frequency}

The zero-frequency case is qualitatively distinct in Theorem~\ref{main 1 second}. When
\begin{equation}
\omega_j=0,
\qquad
j=1,\ldots,N,
\label{eq:num-zero-omega}
\end{equation}
the theorem gives convergence to a stationary configuration without an
additional high-coherence assumption.

We illustrate this regime using a genuinely hybrid system with
\begin{equation}
N=20,
\qquad
n=8.
\label{eq:num-zero-N}
\end{equation}
The positive inertias of the second-order components and the damping
coefficients of the full ensemble span the ranges
\begin{align}
10^{-2}
&\leq m_j \leq 10^{2},
\nonumber\\
10^{-2}
&\leq d_j \leq 10^{2}.
\label{eq:num-zero-multiscale}
\end{align}
The initial phases are distributed over the full circle, and the coupling
strength is $\lambda=1$.

Despite the severe separation of time scales, the trajectory eventually
converges numerically. At $T=3000$, we obtain
\begin{equation}
R_{\rm tail}=1
\label{eq:num-zero-R}
\end{equation}
to the displayed numerical precision, together with
\begin{align}
V_{\rm tail}
&\approx 1.44\times10^{-13},
\nonumber\\
D_{\rm rel}
&\approx 7.46\times10^{-14},
\nonumber\\
D_Z
&\approx 3.66\times10^{-14}.
\label{eq:num-zero-diagnostics}
\end{align}

The transient is strongly multiscale, but all components ultimately approach
a stationary configuration.

\begin{figure}[t]
    \centering
    \includegraphics[width=\textwidth]{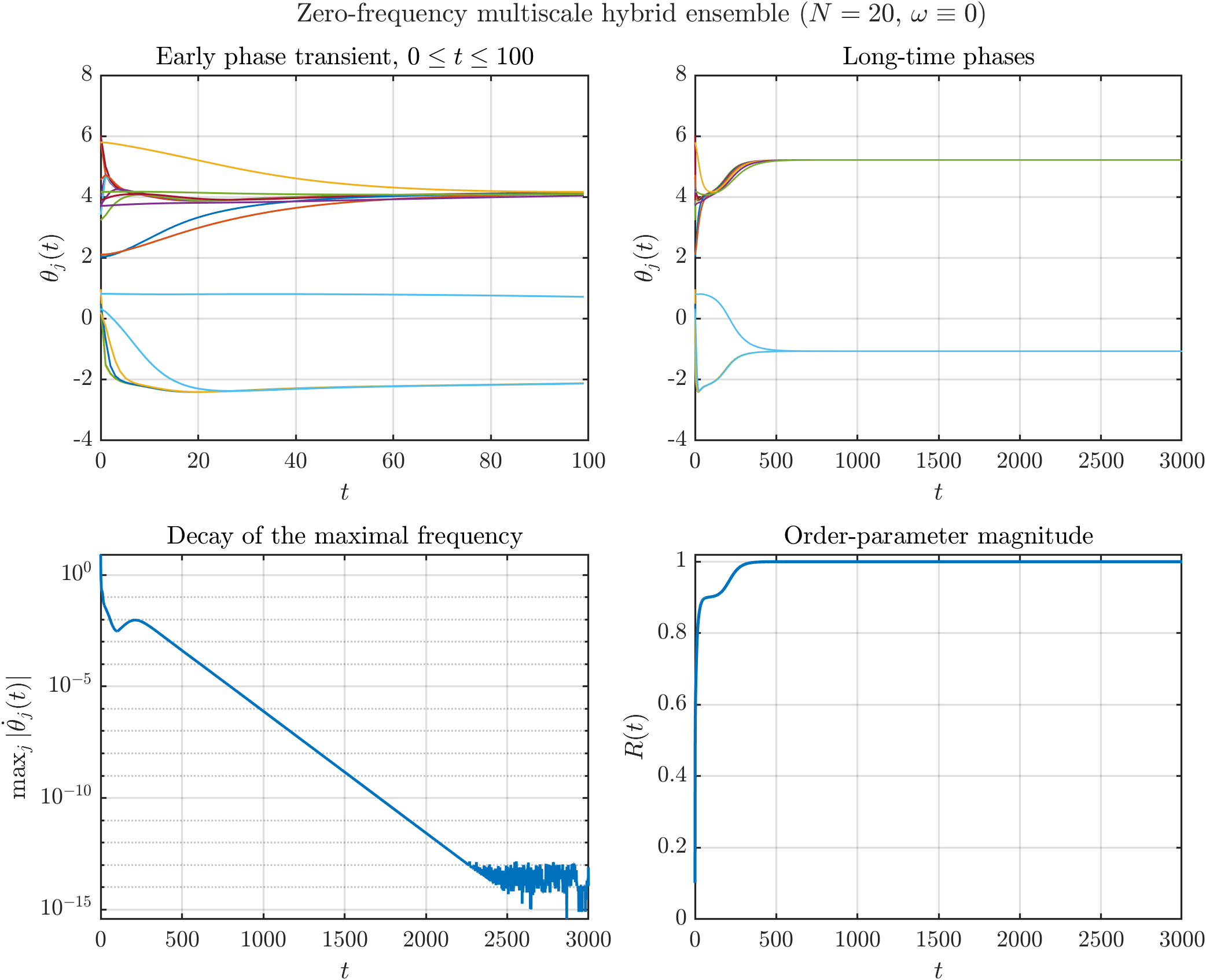}
    \caption{
    Zero-natural-frequency hybrid dynamics with $N=20$ and inertia and
    damping coefficients spanning four orders of magnitude. The early-time
    dynamics exhibit a pronounced multiscale transient, while the long-time
    phase, frequency, and order-parameter diagnostics converge. The figure
    illustrates the separate $\omega=0$ convergence result in the presence of
    degenerate inertia and severe time-scale heterogeneity.
    }
    \label{fig:zero-frequency-multiscale}
\end{figure}

\subsection{Large-scale hybrid dynamics with heterogeneous inertia and damping}
\label{subsec:num-large-N}

Finally, we consider a larger genuinely hybrid ensemble with
\begin{equation}
N=100,
\qquad
n=40.
\label{eq:num-N100}
\end{equation}
For the $60$ second-order oscillators, the positive inertia coefficients are
assigned from the set
\begin{equation}
\left\{
10^{-2},
10^{-1},
1,
10,
10^{2}
\right\},
\label{eq:num-N100-masses}
\end{equation}
and the damping coefficients of the full ensemble are likewise assigned from
\begin{equation}
\left\{
10^{-2},
10^{-1},
1,
10,
10^{2}
\right\}.
\label{eq:num-N100-damping}
\end{equation}

The natural frequencies are obtained from a centered Gaussian realization
with standard deviation $0.3$, for which
\begin{equation}
\omega_M \approx 0.738495,
\label{eq:num-N100-omegaM}
\end{equation}
and the coupling strength is
\begin{equation}
\lambda=7.5.
\label{eq:num-N100-lambda}
\end{equation}

For $N=100$, the geometric high-coherence threshold is
\begin{equation}
1-\frac{2}{N}=0.98.
\label{eq:num-N100-threshold}
\end{equation}
The numerical trajectory nevertheless converges to a highly coherent state,
with
\begin{equation}
R_{\rm tail}\approx0.999291>0.98.
\label{eq:num-N100-R}
\end{equation}
The corresponding terminal diagnostics are
\begin{align}
V_{\rm tail}
&\approx1.29\times10^{-8},
\nonumber\\
D_{\rm rel}
&\approx1.73\times10^{-12},
\nonumber\\
D_Z
&\approx1.16\times10^{-13}.
\label{eq:num-N100-diagnostics}
\end{align}


Under the original integration tolerances, the maximum terminal frequency residual is approximately $1.3 \times 10^{
-8}$
. A numerical refinement study with tighter tolerances and smaller maximum step sizes reduces this residual to approximately $1.4 \times 10^{
-14}$. Thus, the apparent terminal floor is attributable to numerical resolution rather than persistent frequency motion.

\begin{figure}[t]
    \centering
    \includegraphics[width=\textwidth]{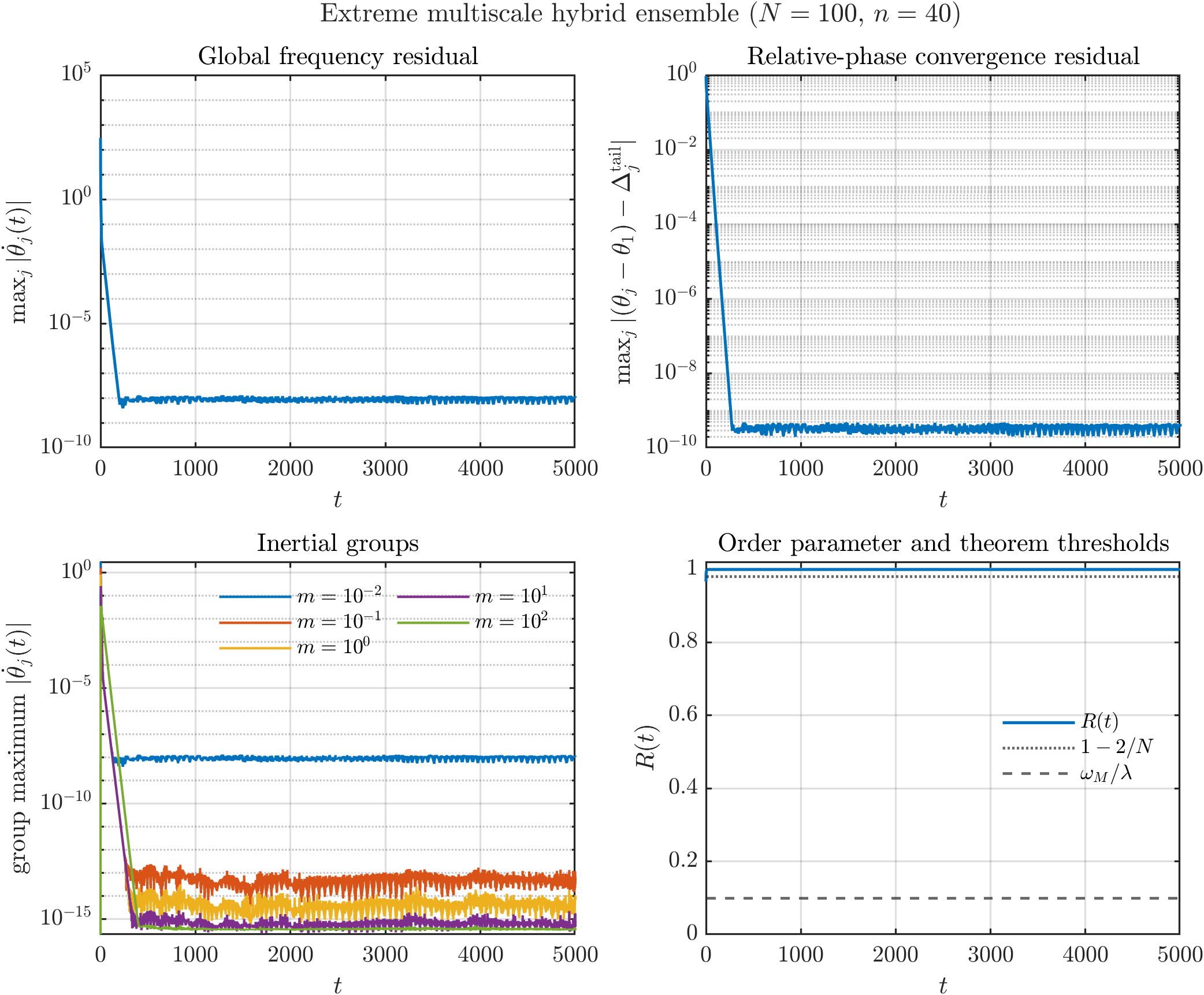}
    \caption{
    Large $N=100$ genuinely hybrid system with both inertia and damping
    spanning $10^{-2}$--$10^{2}$. The grouped frequency residuals reveal
    strongly separated relaxation time scales, while the relative phases and
    order parameter stabilize and $R(t)$ enters the strict high-coherence
    region $R>1-2/N=0.98$. The figure provides a numerical stress test of
    Theorem~\ref{main 1 second} without common-inertia or common-damping assumptions.
    }
    \label{fig:extreme-multiscale}
\end{figure}

The numerical experiments above are intended only to illustrate the
asymptotic relations proved in Theorem~\ref{main 1 second}. In particular, the coupling sweep
shows explicitly that the strict high-coherence condition should not be
interpreted as a necessary synchronization threshold. We make no numerical
or analytical claim here concerning a complete characterization of the
remaining low-coherence regime.
\color{black}

\section*{Acknowledgements} 
{\color{black}The authors thank the associate editor and the two anonymous referees
for their careful reading and constructive comments, which led to
substantial improvements in the manuscript.}
The authors would like to thank Prof. Jared Bronski for insightful discussions during the early stages of this work. T.-Y. Hsiao is grateful to Prof. Shoou-Ren Hsiau for suggestions on the manuscript structure. We are indebted to several friends, Winnie Wang, Ye Zhang, Joy Hsu, Guan-Zhong Chen, and Kui-Yo Chen, for their encouragement and thoughtful input during the preparation of this manuscript. T.-Y. Hsiao is supported by  the European Union  ERC CONSOLIDATOR GRANT 2023 GUnDHam, Project Number: 101124921. Views and opinions expressed are however those of the authors only and do not necessarily reflect those of the European Union or the European Research Council. Neither the European Union nor the granting authority can be held responsible for them.

\section*{Declarations}
The authors declare that they have no conflict of interest. Ethical approval is not applicable.

\section*{Data and Code Availability Statement}
The MATLAB code, complete numerical parameter sets and initial conditions,
and numerical data used to generate the figures and diagnostics in this study
are available from the corresponding author upon reasonable request.

\bibliographystyle{unsrtnat}
\bibliography{sn-bibliography}

@article{strogatz2000kuramoto,
  title={From {K}uramoto to {C}rawford: exploring the onset of synchronization in populations of coupled oscillators},
  author={Strogatz, Steven H},
  journal={Physica D: Nonlinear Phenomena},
  volume={143},
  number={1-4},
  pages={1--20},
  year={2000},
  publisher={Elsevier}
}

@article{ha2016emergence,
  title={Emergence of phase-locked states for the {K}uramoto model in a large coupling regime},
  author={Ha, Seung-Yeal and Kim, Hwa Kil and Ryoo, Seung-Yeon},
  journal={Communications in Mathematical Sciences},
  volume={14},
  number={4},
  pages={1073--1091},
  year={2016},
  publisher={International Press of Boston}
}

@article{dorfler2011critical,
  title={On the critical coupling for {K}uramoto oscillators},
  author={D{\"o}rfler, Florian and Bullo, Francesco},
  journal={SIAM Journal on Applied Dynamical Systems},
  volume={10},
  number={3},
  pages={1070--1099},
  year={2011},
  publisher={SIAM}
}

@article{mehta2015algebraic,
  title={Algebraic geometrization of the {K}uramoto model: Equilibria and stability analysis},
  author={Mehta, Dhagash and Daleo, Noah S and D{\"o}rfler, Florian and Hauenstein, Jonathan D},
  journal={Chaos: An Interdisciplinary Journal of Nonlinear Science},
  volume={25},
  number={5},
  pages={053103},
  year={2015},
  doi={10.1063/1.4919696},
  publisher={AIP Publishing}
}

@book{kuramoto1984chemical,
  title={Chemical Oscillations, Waves, and Turbulence},
  author={Kuramoto, Yoshiki},
  series={Springer Series in Synergetics},
  volume={19},
  publisher={Springer-Verlag},
  address={Berlin, Heidelberg},
  year={1984},
  doi={10.1007/978-3-642-69689-3}
}

@article{wiesenfeld1998frequency,
  title={Frequency locking in {J}osephson arrays: Connection with the {K}uramoto model},
  author={Wiesenfeld, Kurt and Colet, Pere and Strogatz, Steven H},
  journal={Physical Review E},
  volume={57},
  number={2},
  pages={1563},
  year={1998},
  publisher={APS}
}

@article{thumler2023synchrony,
  title={Synchrony for weak coupling in the complexified {K}uramoto model},
  author={Th{\"u}mler, Moritz and Srinivas, Shesha GM and Schr{\"o}der, Malte and Timme, Marc},
  journal={Physical Review Letters},
  volume={130},
  number={18},
  pages={187201},
  year={2023},
  publisher={APS}
}

@article{hsiao2023synchronization,
  title={Synchronization in the quaternionic {K}uramoto model},
  author={Hsiao, Ting-Yang and Lo, Yun-Feng and Wang, Winnie},
  journal={arXiv preprint arXiv:2309.01893},
  year={2023}
}

@article{bronski2020matrix,
  title={A matrix-valued {K}uramoto model},
  author={Bronski, Jared C and Carty, Thomas E and Simpson, Sarah E},
  journal={Journal of Statistical Physics},
  volume={178},
  number={2},
  pages={595--624},
  year={2020},
  publisher={Springer}
}

@article{lohe2009non,
  title={Non-{A}belian {K}uramoto models and synchronization},
  author={Lohe, MA2539317},
  journal={Journal of Physics A: Mathematical and Theoretical},
  volume={42},
  number={39},
  pages={395101},
  year={2009},
  publisher={IOP Publishing}
}

@article{deville2019synchronization,
  title={Synchronization and stability for quantum {K}uramoto},
  author={DeVille, Lee},
  journal={Journal of Statistical Physics},
  volume={174},
  number={1},
  pages={160--187},
  year={2019},
  publisher={Springer}
}

@article{chen2024complete,
  title={Complete and Partial Synchronization of Two-Group and Three-Group {K}uramoto Oscillators},
  author={Chen, Shih-Hsin and Hsia, Chun-Hsiung and Hsiao, Ting-Yang},
  journal={SIAM Journal on Applied Dynamical Systems},
  volume={23},
  number={3},
  pages={1720--1765},
  year={2024},
  publisher={SIAM}
}

@inproceedings{kuramoto1975self,
  title={Self-entrainment of a population of coupled non-linear oscillators},
  author={Kuramoto, Yoshiki},
  booktitle={International Symposium on Mathematical Problems in Theoretical Physics: January 23--29, 1975, Kyoto University, Kyoto/Japan},
  pages={420--422},
  year={1975},
  organization={Springer}
}

@article{hsia2019synchronization,
  title={On the synchronization theory of {K}uramoto oscillators under the effect of inertia},
  author={Hsia, Chun-Hsiung and Jung, Chang-Yeol and Kwon, Bongsuk},
  journal={Journal of Differential Equations},
  volume={267},
  number={2},
  pages={742--775},
  year={2019},
  publisher={Elsevier}
}

@article{coss2018locating,
  title={Locating and counting equilibria of the {K}uramoto model with rank-one coupling},
  author={Coss, Owen and Hauenstein, Jonathan D and Hong, Hoon and Molzahn, Daniel K},
  journal={SIAM Journal on Applied Algebra and Geometry},
  volume={2},
  number={1},
  pages={45--71},
  year={2018},
  publisher={SIAM}
}

@article{baillieul1982geometric,
  title={Geometric critical point analysis of lossless power system models},
  author={Baillieul, John and Byrnes, Christopher},
  journal={IEEE Transactions on Circuits and Systems},
  volume={29},
  number={11},
  pages={724--737},
  year={1982},
  publisher={IEEE}
}

@article{landau1925ungleichungen,
    author = {von Landau, Edmund},
    title = {Einige Ungleichungen F{\"u}r Zweimal Differentiierbare Funktionen},
    journal = {Proceedings of the London Mathematical Society},
    volume = {s2-13},
    number = {1},
    pages = {43-49},
    year = {1914},
    month = {01},
    issn = {0024-6115},
}

@article{xi2017synchronization,
  title={Synchronization of cyclic power grids: equilibria and stability of the synchronous state},
  author={Xi, Kaihua and Dubbeldam, Johan LA and Lin, Hai Xiang},
  journal={Chaos: An Interdisciplinary Journal of Nonlinear Science},
  volume={27},
  number={1},
  pages={013109},
  year={2017},
  doi={10.1063/1.4973770},
  publisher={AIP Publishing}
}

@article{chopra2009exponential,
  title={On exponential synchronization of {K}uramoto oscillators},
  author={Chopra, Nikhil and Spong, Mark W},
  journal={IEEE transactions on Automatic Control},
  volume={54},
  number={2},
  pages={353--357},
  year={2009},
  publisher={IEEE}
}

@article{ha2020asymptotic,
  title={Asymptotic phase-locking dynamics and critical coupling strength for the {K}uramoto model},
  author={Ha, Seung-Yeal and Ryoo, Seung-Yeon},
  journal={Communications in Mathematical Physics},
  volume={377},
  number={2},
  pages={811--857},
  year={2020},
  publisher={Springer}
}

@article{dorfler2013synchronization,
  title={Synchronization in complex oscillator networks and smart grids},
  author={D{\"o}rfler, Florian and Chertkov, Michael and Bullo, Francesco},
  journal={Proceedings of the National Academy of Sciences},
  volume={110},
  number={6},
  pages={2005--2010},
  year={2013},
  publisher={National Academy of Sciences}
}

@article{hsiao2025synchronization,
  title={Synchronization in the complexified {K}uramoto model},
  author={Hsiao, Ting-Yang and Lo, Yun-Feng and Wang, Winnie},
  journal={Nonlinearity},
  volume={39},
  number={1},
  pages={015003},
  year={2026},
  doi={10.1088/1361-6544/ae29d9}
}

@article{dorfler2012synchronization,
  title={Synchronization and transient stability in power networks and nonuniform {K}uramoto oscillators},
  author={D{\"o}rfler, Florian and Bullo, Francesco},
  journal={SIAM Journal on Control and Optimization},
  volume={50},
  number={3},
  pages={1616--1642},
  year={2012},
  publisher={SIAM}
}

@article{hsiaoLoZhu2026equivalence,
  author        = {Hsiao, Ting-Yang and Lo, Yun-Feng and Zhu, Chengbin},
  title         = {Equivalence of Synchronization States in the {Kuramoto} Flow: A Unified Framework},
  journal       = {Journal of Physics: Complexity},
  year          = {2026},
  volume        = {7},
  number        = {2},
  pages         = {025010},
  doi           = {10.1088/2632-072X/ae6be3}
}

@article{acebron2005kuramoto,
  title={The {K}uramoto model: A simple paradigm for synchronization phenomena},
  author={Acebr{\'o}n, Juan A and Bonilla, Luis L and P{\'e}rez Vicente, Conrad J and Ritort, F{\'e}lix and Spigler, Renato},
  journal={Reviews of Modern Physics},
  volume={77},
  number={1},
  pages={137--185},
  year={2005},
  publisher={APS}
}

@article{bronski2012fully,
  title={Fully synchronous solutions and the synchronization phase transition for the finite-{N} {K}uramoto model},
  author={Bronski, Jared C and DeVille, Lee and Park, Moon Jip},
  journal={Chaos: An Interdisciplinary Journal of Nonlinear Science},
  volume={22},
  number={3},
  pages={033133},
  year={2012},
  doi={10.1063/1.4745197},
  publisher={AIP Publishing}
}

@article{bergen2007structure,
  title={A Structure Preserving Model for Power System Stability Analysis},
  author={Bergen, Arthur R. and Hill, David J.},
  journal={IEEE Transactions on Power Apparatus and Systems},
  volume={PAS-100},
  number={1},
  pages={25--35},
  year={1981},
  doi={10.1109/TPAS.1981.316883},
  publisher={IEEE}
}

@article{filatrella2008analysis,
  title={Analysis of a power grid using a {K}uramoto-like model},
  author={Filatrella, Giovanni and Nielsen, Arne Hejde and Pedersen, Niels Falsig},
  journal={The European Physical Journal B},
  volume={61},
  number={4},
  pages={485--491},
  year={2008},
  publisher={Springer}
}

@article{choi2011complete,
  title={Complete synchronization of {K}uramoto oscillators with finite inertia},
  author={Choi, Young-Pil and Ha, Seung-Yeal and Yun, Seok-Bae},
  journal={Physica D: Nonlinear Phenomena},
  volume={240},
  number={1},
  pages={32--44},
  year={2011},
  publisher={Elsevier}
}

@article{choi2014complete,
  title={Complete entrainment of {K}uramoto oscillators with inertia on networks via gradient-like flow},
  author={Choi, Young-Pil and Li, Zhuchun and Ha, Seung-Yeal and Xue, Xiaoping and Yun, Seok-Bae},
  journal={Journal of Differential Equations},
  volume={257},
  number={7},
  pages={2591--2621},
  year={2014},
  publisher={Elsevier}
}

@article{li2014synchronization,
  title={Synchronization and transient stability in power grids based on {{\L}}ojasiewicz inequalities},
  author={Li, Zhuchun and Xue, Xiaoping and Yu, Daren},
  journal={SIAM Journal on Control and Optimization},
  volume={52},
  number={4},
  pages={2482--2511},
  year={2014},
  publisher={SIAM}
}

@article{farkas2016variations,
  title={Variations on {B}arb{\u{a}}lat's lemma},
  author={Farkas, B{\'a}lint and Wegner, Sven-Ake},
  journal={The American Mathematical Monthly},
  volume={123},
  number={8},
  pages={825--830},
  year={2016},
  publisher={Taylor \& Francis}
}

@book{kolmogorov1949inequalities,
  title={On Inequalities Between the Upper Bounds of the Successive Derivatives of an Arbitrary Function on an Infinite Interval},
  author={Kolmogorov, Andrey Nikolaevich},
  publisher={American Mathematical Society},
  address={New York},
  year={1949}
}

@article{cho2025inertia,
  title={Inertia perturbation theory for the inertial {K}uramoto model},
  author={Cho, Hangjun and Dong, Jiu-Gang and Ha, Seung-Yeal and Ryoo, Seung-Yeon},
  journal={SIAM Journal on Mathematical Analysis},
  volume={58},
  number={4},
  pages={3669--3719},
  year={2026},
  doi={10.1137/25M1788567}
}

@article{yeung1999time,
  title={Time delay in the {K}uramoto model of coupled oscillators},
  author={Yeung, MK Stephen and Strogatz, Steven H},
  journal={Physical Review Letters},
  volume={82},
  number={3},
  pages={648},
  year={1999},
  publisher={APS}
}

@article{turner1998five,
  title={Five parametric resonances in a microelectromechanical system},
  author={Turner, Kimberly L and Miller, Scott A and Hartwell, Peter G and MacDonald, Noel C and Strogatz, Steven H and Adams, Scott G},
  journal={Nature},
  volume={396},
  number={6707},
  pages={149--152},
  year={1998},
  publisher={Nature Publishing Group UK London}
}

@article{hardin1990feedback,
  title={Feedback of the Drosophila period gene product on circadian cycling of its messenger {RNA} levels},
  author={Hardin, Paul E and Hall, Jeffrey C and Rosbash, Michael},
  journal={Nature},
  volume={343},
  number={6258},
  pages={536--540},
  year={1990},
  publisher={Nature Publishing Group UK London}
}

@article{bargiello1984restoration,
  title={Restoration of circadian behavioural rhythms by gene transfer in {D}rosophila},
  author={Bargiello, Thaddeus A and Jackson, F Rob and Young, Michael W},
  journal={Nature},
  volume={312},
  number={5996},
  pages={752--754},
  year={1984},
  publisher={Nature Publishing Group UK London}
}

@article{mirollo1990synchronization,
  title={Synchronization of pulse-coupled biological oscillators},
  author={Mirollo, Renato E and Strogatz, Steven H},
  journal={SIAM Journal on Applied Mathematics},
  volume={50},
  number={6},
  pages={1645--1662},
  year={1990},
  publisher={SIAM}
}

@article{strogatz2012sync,
  title={Sync: {H}ow order emerges from chaos in the universe, nature, and daily life},
  author={Strogatz, Steven H},
  journal={Grand Central Publishing},
  year={2012}
}

@article{liXue2019gradient,
  author  = {Li, Zhuchun and Xue, Xiaoping},
  title   = {Convergence of analytic gradient-type systems with periodicity and its applications in {Kuramoto} models},
  journal = {Applied Mathematics Letters},
  volume  = {90},
  pages   = {194--201},
  year    = {2019},
  doi     = {10.1016/j.aml.2018.10.015}
}

@article{hirschSmithZhao2001,
  author  = {Hirsch, Morris W. and Smith, Hal L. and Zhao, Xiao-Qiang},
  title   = {Chain transitivity, attractivity, and strong repellors for semidynamical systems},
  journal = {Journal of Dynamics and Differential Equations},
  volume  = {13},
  number  = {1},
  pages   = {107--131},
  year    = {2001},
  doi     = {10.1023/A:1009044515567}
}

@article{hurley1995chain,
  author  = {Hurley, Mike},
  title   = {Chain recurrence, semiflows, and gradients},
  journal = {Journal of Dynamics and Differential Equations},
  volume  = {7},
  number  = {3},
  pages   = {437--456},
  year    = {1995},
  doi     = {10.1007/BF02219371}
}

@article{daiFiedlerLopezNieto2025,
  author        = {Dai, Jia-Yuan and Fiedler, Bernold and
                   L{\'o}pez-Nieto, Alejandro},
  title         = {Transient Rebellions in the {Kuramoto} Oscillator:
                   {Morse--Smale} Structural Stability and Connection
                   Graphs of Finite 2-Shift Type},
  journal={arXiv preprint arXiv:2512.02937},
  year={2025}
}
\end{document}